\input amstex
\documentstyle{amsppt}
\input epsf
\input psfrag

\define\NN{\Bbb N}
\define\QQ{\Bbb Q}
\define\RR{\Bbb R}
\define\CC{\Bbb C}
\define\SS{\Bbb S}

\define\ZZ{\Bbb Z}

\define\PP{\Bbb P}
\define\cat{\mbox{cat}}

\hyphenation{Liou-ville}
\renewcommand{\AA}{\Bbb A}

\NoBlackBoxes

\topmatter
\title Homogeneous Systems and Euclidean Topology
\endtitle
\rightheadtext {}
\author Jon A. Sjogren
\endauthor
\affil Towson University
\endaffil
\address Towson, Maryland
\endaddress
\date 24 July 2017
\enddate
\endtopmatter
\document

\abovedisplayskip=13pt
\belowdisplayskip=13pt

\head Introductory Remarks on the Topology of Euclidean Space\endhead

Fundamental facts that are characteristic of finite dimensional Euclidean space $\RR^n$, a real vector space endowed with the Pythagorean metric include Invariance of Domain (that a locally injective mapping is open) and  the Jordan-Brouwer theorem (that a topological $\SS^{n-1}$ as a subspace of standard $\RR^n$, when removed leaves one bounded and one unbounded component).

These theorems have been approached from  several points of view. Certainly, Brouwer's Fixed-Point theorem, with generalizations, is a powerful implement toward these important results. In addition, there exist now several proofs that use only elementary calculus and are easy to comprehend. The theorem known by names of H. Poincar\'e and C. Miranda (but understood earlier by Hadamard, Kronecker and others), is ``equivalent'' to BFPT and sometimes  makes a more direct application to the problem at hand.

Also renowned is the  Borsuk-Ulam Theorem in $n$ dimensions. This theorem directly implies BFPT, so it may well the correct tool to use. This utility has been observed more often in texts on non-linear analysis, see \cite{Deimling} than those on topology. So we aim at a  suitable proof of the Borsuk-Ulam (B-U) Theorem. The various versions of B-U will not  formally be listed; they can be found, together with the Lusternik-Schnirel'mann covering theorem in the book of \cite{Matou{\v s}ek} and the notes of \cite{Suciu}. We wish to avoid most of the proofs commonly cited, that require high-powered theory, complicated constructors, or  subtle concepts that are extraneous to the  problem at hand.

The pathway we choose starts with transforming the problem from one of ``continuous mapping'' to  one of ``solve a collection of homogeneous multi-nomials'' by means of the Weierstra{\ss} Approximation Theorem. The latter result is quite effective as seen from an analytic solution to the Heat equation, or one of the  formulas that yield the multinomial coefficients, such as the expressions due to Bernstein or Landau, see \cite{Sjogren, Iterated}.

It turns out that we have arrived at a purely algebraic problem exposited by A. Pfister. The result actually has meaning for any real-closed field $R$ that is ground field to a vector space, not only for the standard reals $\RR$. For those topological analysts whose facility in the homological theory of commutative rings may not rise to the level achieved by Prof\. Pfister,  there is a way to simplify the proof in the  case of the standard reals, as noted in \cite{Lang, Places} and by others. The Annals of Mathematics paper of S. Lang does not directly refer to the  B-U Theorem however.

We express the Borsuk-Ulam Theorem in a minimalist form: any {\it odd} (antipode-preserving) mapping $\varphi: \SS^n \to \SS^n$ is {\it essential} (meaning not nil-homotopic, not contractible within the image space $\SS^n$). Of course the version stating that the Brouwer degree $\deg\varphi$ is an odd integer, is sharper. A statement, only {\it apparently} more general, is that a ``$\ZZ_2$-equivariant'' self-map of the sphere is essential.

Our version of B-U allows one immediately to proceed to the Invariance of Domain Theorem without using any numerical invariants. In particular we avoid formulas for the computation of mapping degree. Also we will not have used the Leray Product Formula, Borsuk's separation thesis, various simplicial approximations, or characterization of  connected components, amongst other subtle concepts of general topology. We say ``subtle'', noting in particular that classical Domain Invariance from \cite{Hurewicz \& Wallman} was  covered in \cite{Dugundji}, but has not further been explained in recent texts, except by rote. Prof. T. Tao gives a concise proof in his on-line journal \cite{Tao, blog},  related to work by W. Kulpa, using  only metric topology.

\bigskip
\head Borsuk-Ulam via Projective Varieties\endhead

Our main supporting result is a polynomial verion of B-U in the spirit of \cite{Knebusch}, \cite{Arason} and \cite{Pfister}. This result applies to any real-closed ground field $R$, not only the standard reals. From there to reach the usual B-U theorem, then to Invariance of Domain, we need basic observations about $\RR$ of an analytic nature. This is analogous to the saying, ``the Fundamental Theorem of Algebra (that a real polynomial of degree $\geq 3$ is reducible) cannot be proven using algebra only''.  In fact, given Brouwer's Invariance of Domain, a topological proof of ``FTA'' is readily derived, \cite{Sjogren, Domain}, in the {\it real} form as stated, not mentioning complex numbers.

In other words, a completeness property of $\RR$ such as the Bolzano theorem on the least upper bound is required. Furthermore, one must employ ``compactness'' in the sense that ``the space of lines (real projective space) is compact'' in finite dimension. This tells the Analyst that an accumulation point of a subspace lies in the projective space. Hence we may find a point $y \in \SS^n$ that maps to $\vec{0} \in \RR^n$ by $\varphi: \SS^n \to \RR^n$, a given continuous mapping.

In the standard case of ground field $R$, Prof. Lang could simplify his quasi-real B\'ezout theorem for polynomials (``multi-nomials'') of {\it odd} total degree. The connection to the Borsuk-Ulam question was not understood until later.

So let us state a ``multi-nomial version'' of the B-U Theorem, and indicate a proof by algebraic methods.  Then it is not surprising that the Weierstra{\ss} approximation leads to the full ``continuous'' B-U result.

\bigskip
\noindent
{\bf Theorem 1}\quad  Given a quantity $n$ of polynomials over the field $\RR$ in $n+1$ variables $q_1(x_1,\dotsc, x_{n+1}),\dotsc, q_n (x_1,\dotsc, x_{n+1})$, which are all ``odd'', namely
$$q_j\left(-x_1, -x_2, \dotsc, -x_{n+1}\right) = -q_j \left(x_1,\dotsc, x_{n+1}\right)$$
for $j=1,\dotsc, n$. Then there exists a {\it ray} consisting of all vectors $\lambda (b_1,\dotsc, b_{n+1})$ where $\lambda > 0$ and $\vec{b}$ is not the zero vector, such that $q_j (\lambda \vec{b}) = 0$ for each $j =1,\dotsc, n$.

\bigskip
We will see that in Theorem 1, the standard reals $\RR$ can be replaced by any other real-closed field $R$. We defer the proof until Theorem 3 has been stated, which is actually the principal result of  the Section. It seems that the general B-U theorem for $\alpha$ continuous mapping is {\it not true} for any real-closed field $R$ other than $\RR$.

\bigskip
\noindent
{\bf Theorem 2 (Borsuk-Ulam)}\quad
Given $f: \SS^n \to \RR^n$, an odd continuous mapping, so that for $X \in \SS^n$ there holds $f(-x) = -f(x)$, then there exists $y \in \SS^n$ such that $f(y) = \vec{0} \in \RR^n$.

\bigskip
\noindent
{\bf Remark}\quad Another formulation is that {\it any} continuous  $g: \SS^n \to \RR^n$ yields up some $y \in \SS^n$ satisfying $g(y) = g(-y) \in \RR^n$.

\bigskip
\noindent
{\it Sketch of proof}\quad Let $f_i: \SS^n \to \RR$, $i =1,\dotsc, n$ be the coordinates of $f$. Then taking $\{p_i(x)\}$ to be $\epsilon$-approximations to $\{f_i(x)\}$, where $p_i(x)$ is the restriction to $\SS^n$ of a real multi-nomial $p_i(x_1,\dotsc, x_{n+1})$, we may actually replace $p_i(x)$ by $q_i(x) = \frac{1}{2} \left[p_i(x)-p_i(-x)\right]$ and obtain for all $1 \leq i \leq n$
$$\left|f_i(x) -q_i(x)\right| < \epsilon\quad\mbox{for}\quad x \in\SS^n.$$
This inequality holds since $f_i$ is odd, and we also know that $q_i(x)$ is odd from its definition. Next, if $f_i$ on $\SS^n$ is bounded away from $\phi$ by $\delta >0$, then all $\{q_j\}$ are bounded away from $\phi$ in modulus by $\delta-\epsilon >0$, where we chose $\epsilon >0$ small enough. By continuity of $\{q_j\}$ and  compactness of $\SS^n$, one may infer that the $\{f_j\}$ have no common zero (as a ray), which contradicts Theorem 1.\hfill $\blacksquare$

\bigskip
This proof uses the well-known (to analysts) ``compactness argument'' whereby a sequence of values in a compact space gives rise to a ``convergent sub-sequence'' or equivalently an ``accumulation point''. We need to use the compactness argument again in this section.

We now state the form of B\'ezout's theorem ``over a real-closed field'' that is required. We use the standard real numbers $\RR$ as our prototype or main exemplar of a real-closed field.

\bigskip
\noindent
{\bf Theorem 3}\quad For $n \geq 1$, let $f_1,\dotsc, f_n \in \RR \left[x_1,\dotsc, x_{n+1}\right]$ be homogeneous multi-nomials (forms) of respective degrees $d_1, \dotsc, d_n$, with each $d_i$ an {\it odd} natural number. Then there exists a non-zero real solution vector $\vec{a} = (a_1,\dotsc, a_{n+1}) \in \RR^{n+1}$, that is, satisfying $f_j(a_1,\dotsc, a_{n+1}) = 0$ for $j=1,\dotsc, n$. In fact, $\vec{a}$ generates a solution ray $\{\lambda \vec{a}\}$, $\lambda \ne 0$, $\lambda$ real.

\bigskip
\noindent
{\bf Remark}\quad Several components of a proof are indicated, which may be selected and assembled according to  the taste of the reader. The proof should be ``algebraic enough'' still to hold for other real-closed fiels.

\bigskip
The ``simplest'' proof is perhaps constituted by the observation that Theorem 3 is exactly the Theorem 1 given on page 239 of \cite{Shafarevich}.

Thus the reader who accepts certain results ``modulo the algebra'' now has  the Borsuk-Ulam theorem fully in hand (once the derivation of our Theorem 1 is completed as a Corollary). The treatment in \cite{Shafarevich} is straightforward based on the theory of algebraic divisors. Nevertheless, we proceed to redo parts of this work based on the concept of Resultant Systems (\cite{Macaulay}, \cite {Kapferer}, \cite{vd Waerden 1927}, \cite{Behrend}), which embodies Algebraic Geometry of a generation or two prior to Basic Algebraic Geometry, Vol. I. In volume II the same learned author Prof. Shararevich treats the  contemporaneous theory of schemes developed by Serre-Grothendieck. In the continuation, which is largely based on B.L. van der Waerden's foundational articles and chapters, we intend for definitions to be reasonably concrete. For example, ``multiplicity of a solution'' should be calculable from  Polynomial Ideal Theory.

Now recall that we wished to establish B-U theorem at least for multi-nomial functions.

\bigskip
\noindent
{\it Proof of Theorem 1}\quad See \cite{Pfister}. We have a quantity $n$ of polynomials $\{q_j\}$ which are {\it odd} as functions in their $n+1$ arguments $x_1,\dotsc, x_{n+1}$, but we may homogenize the $q_j$ by throwing in an additional variable to achieve the required total degree.

For example, $q(x_1, x_2, x_3) = 2x_1-x_2x_3^2+x_1^3x_2x_3-3x_1x_3^2+x_2^2x_3$ satisfies $q(-x_1, -x_2,-x_3) = -q(x_1,x_2,x_3)$. Note that the degree of each term is odd, so the needed power of $x_0$ is always {\it even}. Take
$$\tilde q = (x_0, x_1, x_2, x_3) = 2x_0^4x_1-
x_0^2x_2x_3^2+
x_1^3x_2x_3-
3x_0^2x_1x_3^2-
x_0^2x_2^2x_3.$$
It is not difficult to show that the above observation on degrees holds in general. Now for each $j = 1,\dotsc, n$, replace any factor $x_0^2$ by $x_1^2 +\cdots+x_{n+1}^2$ in $q_j$. Doing so yields a quantity $n$ of  odd-degree homogeneous polynomials (or multi-nomials) $\hat{q}_j(x_1,\dotsc, x_{n+1})$ which by Theorem 3 above possess a common solution valid on a ray in $\RR^{n+1}$ that is generated by a non-zero real vector $(a_1,\dotsc, a_{n+1})$. By homogenity of the $\hat{q}_j$, we may choose the solution vector $\vec{b} = \frac{\vec{a}}{\|\vec{a}\|} \in \SS^n$.

Also $-\vec{b}$ is an acceptable solution. Either can be taken as the point on the $n$-sphere (or on $\RR P^n$) sought by the Borsuk-Ulam theorem (expressed also in Theorem 2).\hfill $\blacksquare$

\bigskip
\head Concerning the ``homotopy'' interpretation of B-U Theorem\endhead

Strong versions of the theorem exist, in the form of ``an antipode-preserving mapping $g: \SS^n \to \SS^n$ has odd Brouwer degree''.

\bigskip
\noindent
{\bf Homotopy Borsuk-Ulam Theorem}\quad Such an {\it odd} mapping (commuting with  the canonical involution of $\SS^n$) is {\it essential}. That is, $g$ is not contractible to a point in the image sphere $\SS^n_w$. The conclusion once again is that $g$ is {\it not} homotopic within $\SS^n_w$ to any constant mapping.

\bigskip
\noindent
{\it Proof}\quad We deduce this from Theorem 2. Also the result {\it implies} Theorem 2 directly, \cite{Matou{\v s}ek}. For $g$ to be inessential or nil-homotopic means that there is an extension $\tilde{g}: B^{n+1}\to \SS^n$ of $g$ whose domain is the Euclidean ball $B^{n+1}$ with boundary $\SS^n$.  That is, $\tilde{g}$ restricted to $\partial B^{n+1}$ is just $g$, see \cite{Dugundji}. Next we may define the projection $\pi: \SS^{n+1}_+ \to B^{n+1}$ from the ``upper hemisphere'' of $\partial B^{n+2}$ by means of $\pi(x_1,\dotsc, x_{n+2}) = (x_1, \dotsc, x_{n+1})$ where $x_{n+2} > 0$ and $\sum_{i=1}^{n+2} x_i^2 = 1$. Thus we have a continuous  mapping $f: \tilde{g} \circ \pi: \SS^{n+1} \to \SS^n$ and similarly $f: \SS^{n+1} \to \SS^n$ on the lower hemisphere, defined by $f(x) = -\tilde{g} \circ \pi(-x)$. Since $\tilde{g}_{|\SS^n}$ is antipode-preserving (odd), the mapping $f: \SS^n \to \SS^n$ is well-defined, continuous and antipode-preserving, hence it is  also such a mapping $\SS^{n+1} \to \RR^{n+1}$ not meeting the origin, which violates Theorem 2.\hfill $\blacksquare$

\bigskip
For future use, we note a simple 

\noindent
{\bf Homotopy Fact}: suppose that for $g, h:\,\SS^{n-1}\to\SS^{n-1}$, $g \sim h$ (considered as mappings to $\RR^n$) by a homotopy $H: \SS^{n-1} \times I\to \RR^n$. Then if $H(s,t)$ never attains $\vec{0} \in \RR^n$, where $s \in \SS^n$, $t \in [0,1]$, then $h$ is homotopic to $g$ considered as mappings to $\SS^{n-1}_w$.

\bigskip
\noindent
{\it Proof}\quad If $H$ exists, it may be modified by pushing away from $\vec{0}$ and $\infty$ so that all of its values lie on $\SS^{n-1}$. Thus we have a homotopy $\tilde{H}: g \sim h$ within $\SS^{n-1}$. In particular $g$ is essential if and only if $h$ is essential.\hfill $\blacksquare$

\noindent
One {\it consequence} of this {\bf Fact} is that a mapping $g: B_z \to B_w$ of one ball to another ball, which restricts to and {\it essential} map $\partial g : \partial B_z \to \partial B_w$ must itself be surjective onto $B_w$.
\bigskip
\head Classical Domain Invariance\endhead

Background for Brouwer's Invariance of Domain can be found in \cite{Dugundji}, \cite{Deimling} and \cite{Tao, blog}. This famous theorem on the topology of Euclidean space, from around 1910, can be stated:

\bigskip
\noindent
{\bf Theorem IVD1}\quad Let $\Omega \subset \RR^n$ be an open set. Then any (continuous) mapping $h: \Omega \to \RR^n$ that is locally one-to-one, is an open mapping.

\bigskip
By way of explanation of the terminology, we quote an equivalent but more concrete statement.

\bigskip
\noindent
{\bf Proposition IVD2}\quad  Suppose  $g: B^n_z \to B^n_w$ is a one-to-one mapping with $g(\vec{0}_z) = \vec{0}_w$. Then there exists an open subset $U \subset g (B_z)$ with $\vec{0} \in U$. Here $B_z$, $B_w$ are the open unit balls at the Origin, distinguishing ``domain'' from ``range''.
Taking a ball of smaller radius, we could regard $g$ as defined and continuous on $\overline{B_z(1)}$, the closed unit ball.
\bigskip
\noindent
{\bf Proof} \quad Now consider the homotopy
$$H: \overline{B}_z \times I \to \RR^n_w$$
defined by
$$H(x, t) = g \left(\frac{x}{1+t}\right) - g\left(\frac{-tx}{1+t}\right).$$
For all $0 \leq t \leq 1$, $H$ maps the Origin $\vec{0}_z$ to the Origin $\vec{0}_w$. Also $\mbox{Im}(H) \subset \overline{B}_w(1)$, though to avoid a calculation, $H$ could be scaled radially so that  its image fits into $\overline{B}_w(1)$. An important fact is that for $x \in \partial B(\rho)$, $\rho > 0$, $H(x,t)$ is  never $\vec{0}_w$: the homotopy restricted to any sphere of radius $\leq 1$ cannot cross the origin. This follows from the assumption of injectivity for $g$. Thus on each ``central sphere'' $\SS^{n-1}_p = \partial B^n(\rho)$, the mapping $g$ is homotopic to
$$\phi(x) = H(x,1) = g\left(\frac{x}{2}\right)-g\left(\frac{-x}{2}\right)$$
by the  restriction of $y = H(x,t)$, where as $t$ varies, $y$ never crosses the Origin.

We note that $\phi$ on every central sphere $\SS^{n-1}_{\rho}$ is {\it odd} $(\ZZ_{2}$-equivariant or antipode-preserving). By compactness of $\SS^{n-1}(1)$, $g$ and $\phi$ attain their infinum in norm $\|g(x)\|$ and
$\|\phi(x)\|$, $x \in \SS^{n-1}_z (1)$. Choose a radius $\sigma > 0$ smaller than both of these positive infina. Next we consider a deformation retraction $G: B^n_w \times I \to B^n_w (\sigma)$ given by
$$ G(y,t) =    \left\{\matrix
\bigl[\,\sigma t + \|y\| (1-t)\,\bigr]\dfrac{y}{ \|y\|}, & \mbox{for $\|y\| \geq \sigma,$} \\
 y,  & \mbox{for $\|y\| < \sigma$.}
\endmatrix\right.$$
One notes that $G$ is ``piece-wise linear'' and not generally smooth on $\SS_a(\sigma)$. In the following Figures we suppress the dimensions of the Spheres and other spaces that are depicted.

\vfill
\centerline{\epsfbox{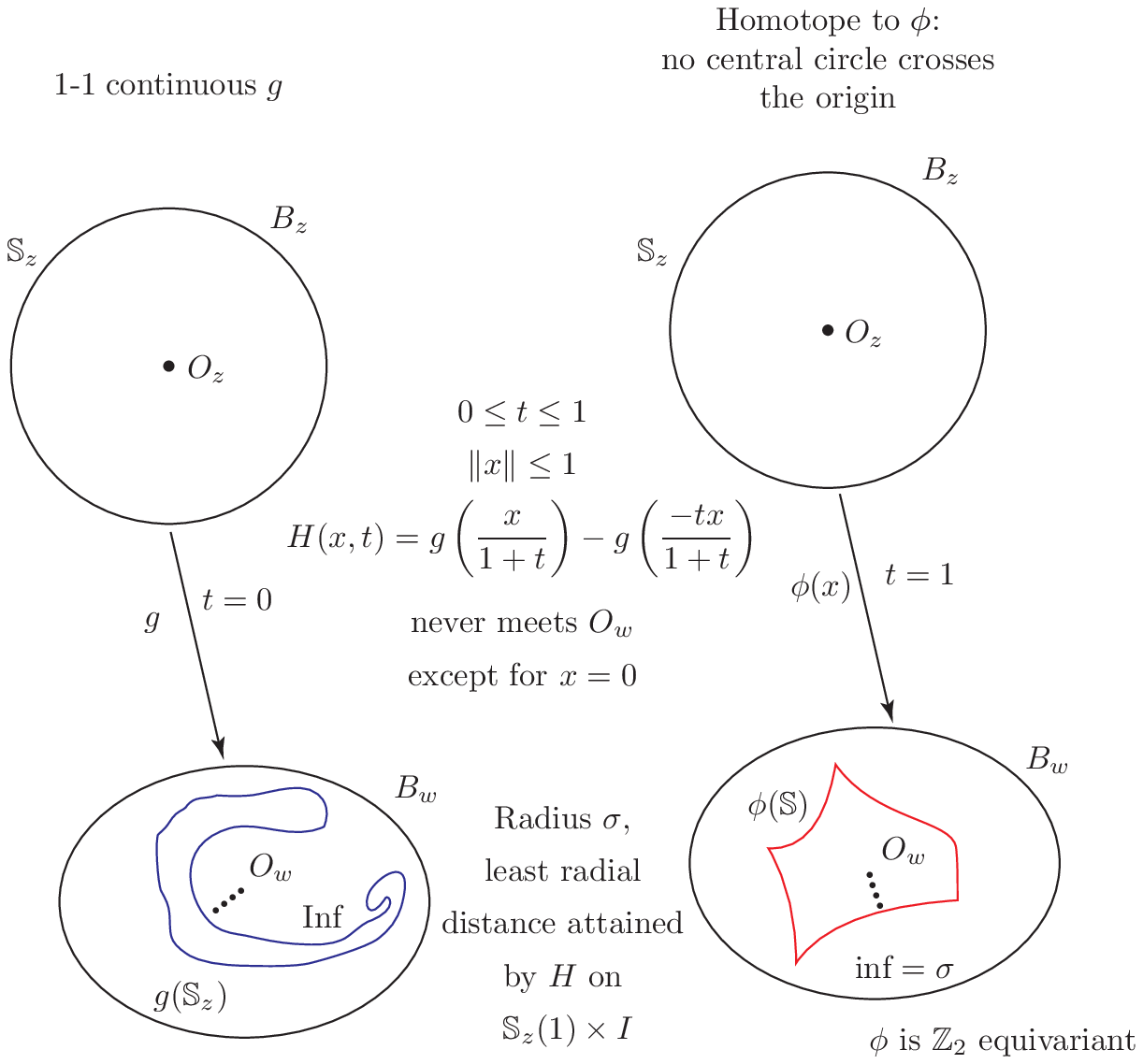}}

\medskip

\centerline{\bf Figure 1}

\bigskip
For each $t \in I$, the radial ray containing $y$ is kept invariant (in terms of its $z$- and $w$-coordinates). The mapping $G$ is a homotopy between the ``identity'': $B_z^n (1) \to B^n_w(1)$ and the ``radial retraction'': $B^n_z(1)\to B^n_w(\sigma)$ that keeps the smaller ball point-wise fixed.

Now define $L_g(x) = G(y, 1) \circ g(x)$ and $L_{\phi} (x) = G(y,1) \circ \phi(x)$, both of which map $\overline{B_z^n}(1)$ to $\overline{B^n_w}(\sigma)$. Furthermore both $L_g$ and $L_{\phi}$, restricted to $\partial \overline´{B^n_z}(1)$, have\linebreak

 \newpage

 \vglue-1cm
 \centerline{\epsfbox{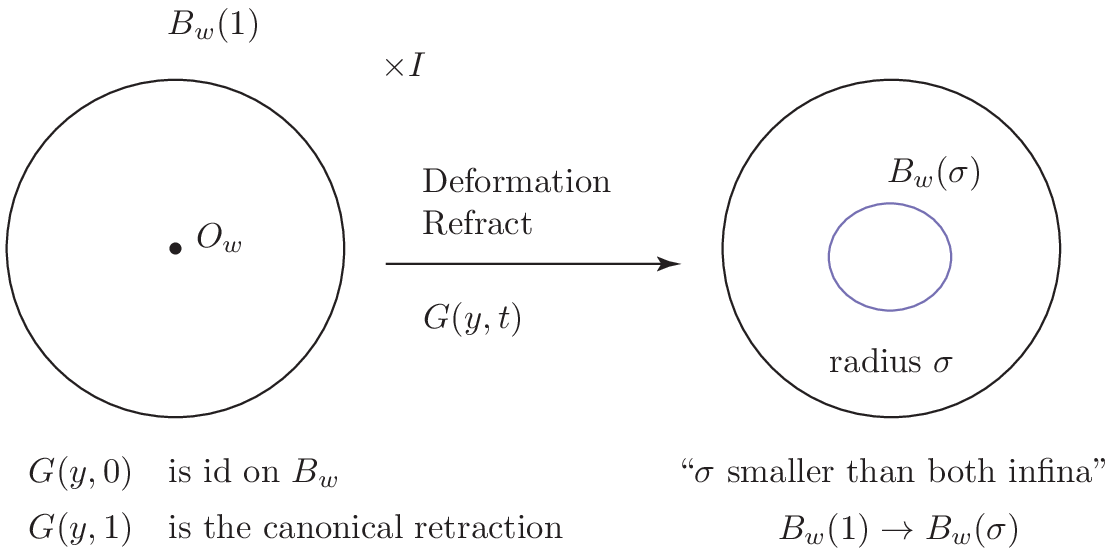}}

\centerline{\bf Figure 2}

\vfill
\centerline{\epsfbox{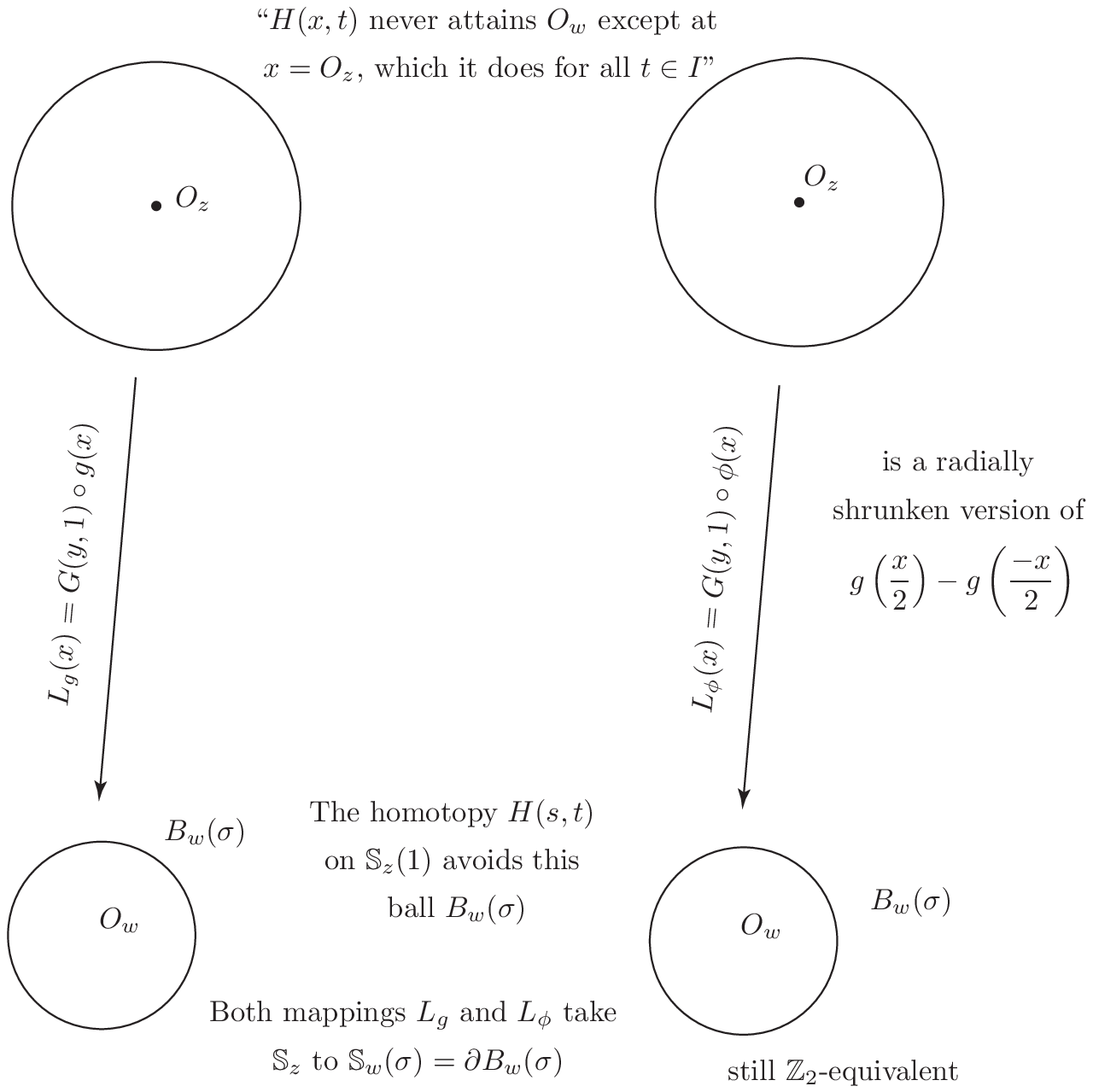}}

\centerline{\bf Figure 3}

 \newpage
 \noindent
 image contained in $\partial \overline{B_w^n}(\sigma)$, the ``small image sphere''. We observe that $L_{\phi}$ restricted to the ``big $z$-sphere'' $\SS^{n-1}_z(1)$ is actually antipode-preserving (also if restricted to other central spheres). Hence by the  homotopy form of the B-U theorem above,  we conclude that $L_{\phi}: \SS^{n-1}_z(1) \to \SS^{n-1}_w(\sigma)$ is an essential mapping.

Also $L_g: \SS^{n-1}_z(1) \to \SS^n_w(\sigma)$, though not an injective mapping, is homotopic to $L_{\phi}$ within this space of mappings, since the homotopy $H(x,t)$ followed by $G(y,s)$ avoids the $w$-Origin, so these homotopies may be projected radially to the $w$-sphere of radius $\sigma$. It follows that $L_g$ restricted to $\SS^{n-1}_z(1)$ is also essential and by the Homotopy Fact above,
$L_g: \overline{B_z^n}(1) \to \overline{B^n_w}(\sigma)$ is a surjection.

A given $b \in B_w(\sigma)$ is therefore in the image of $L_g = G \circ g$, but it is not moved under $G(\cdot\,, t)$. Hence $b = g(a)$ for some $a \in B_z(1)$. Since $b$ was chosen arbitrarily, we have found an open neighborhood $B_w(\sigma)$ of $\vec{0}_w$ in the image, confirming that $g$ must be an open mapping. \hfill $\blacksquare$

\bigskip
\head B\'ezout's Theorem and Solution Multiplicity\endhead

Theorem 3 above follows from Theorem 4, which allows a more general ``ground field'', see  \cite{vd Waerden, Algebra II}, section 83.

\bigskip
\noindent
{\bf Theorem 4 (B\'ezout)}\quad If a system $F_1,\dotsc, F_n$ of homogeneous equations $(F_j=0)$, in $n+1$ variables  $x_1,\dotsc, x_n$, with $x_j \in R(\sqrt{-1})$, with coefficients in the real-closed field, has only finitely many distinct solutions $(x_i)\ne \vec{0}$, then there holds a formula for their multiplicity. Consider as before the solutions generating lines (or {\it rays}) over $C = R(\sqrt{-1})$. Defining
$$\Delta = \prod_{i=1}^n \deg F_j,$$
we obtain
$$\Delta = \sum_{P} \mbox{mult} (P),\tag{B}$$
where $\{P\}$ runs though the distinct solution rays and $\mbox{mult}(P)$ is the {\it multiplicity} of the solution to the given and algebraic definition.

\bigskip
Finally, we are looking for Theorem 3 as a corollary. We may write Theorem 3  again as: 

\noindent
{\bf Theorem 5} \quad With the hypotheses of Theorem 4, given the homogeneous system $F_1,\dotsc, F_n$ with coefficients in the real-closed $R$, suppose that each degree $(F_j) = d_j$ is odd, $j = 1,\dotsc, n$, then we conclude that there exists a solution $(\xi_0: \xi_1: \cdots : \xi_n)$ defining a {\it ray}, with all $\xi_j \in R$.
\bigskip
\noindent
To finish a proof of Theorem 3, we specialize $R$ in Theorem 5 to the ``standard'' real numbers $\RR$. For our purposes, we also need only consider the case (as in \cite{vd Waerden, Algebra II} p. 16, where there exist only finitely many solution (rays). Furthermore, for the application to the B-U theorem, we may assume that  the {\it coefficients} of the equation system are {\it transcendental}, and algebraically independent over the rationals $\QQ$.\hfill $\blacksquare$

A result similar to Theorem 5 has been considered from several points of view, as is seen in the section below.

\bigskip
\head Algebra and Topology in Theorem Five\endhead

We mentioned that on Chapter III of \cite{Shafarevich}, Book 1, the theory of the {\it divisor class group} of a variety is applied to prove B\'ezout's theorem in the form we need, our Theorem 3 or 4. As a matter of fact, this {\it algebraic} method uses a  {\it general position} argument concerning the equations $F_1, \dotsc, F_n$ of our system (which is fulfilled if  then coefficients are algebraically independent or {\it generic}). For the standard real numbers $\RR$ as coefficients,  the usual limiting inference (by compactness $\RR P^n$) gives Theorem 5 more generally. This discussion shows in a rough manner the trade-off between the power of  using the metric on $\RR$, and achieving the Theorem for an arbitrary system (not necessarily generic).

Using divisors on a variety was originally beyond our scope, so we examine proofs that use algebraic geometry of a nature even more elementary. Now a rather pure form of Theorem 5, wholly algebraic in statement and proof, is given in the book of \cite{ Pfister}, p\. 57. The author's remarks point toward an interpretation into geometry of his module-theoretic argument (valid for any real-closed $R$). It is argued that, the greater degree to which the proof is intuited geometrically,  the less it is convincing in its rigor.

Our point of view is that by throwing in a bit of the order or the topology of $\RR$, we obtain a proof of Theorem 5 that is {\it predominantly} algebraic but uses commonly known facts. On the other hand, for $\RR$, work of Borsuk and Hopf from the 1930s on the B-U theorem itself, leads to a purely topological proof (with almost nothing about polynomials). Readers are invited to revisit this part of the history, \cite{Hopf}, where the demonstrations may not be obvious to the  contemporary scholar. By means of the modern machinery of algebraic topology, such proofs can be down-sized; we indicate the section ahead covering the earlier work of L. Lusternik and L\. Schnirel'mann.

\looseness=-1
An early algebraic proof of Theorem 5 is reputed to be that of \cite{Behrend}. Here the result is stated for coefficients in $\RR$, which is our case of interest. The author is looking at our system $F_1,\dotsc, F_n$ where the latter are dependent on several rows of indeterminates, not only $x_0, \dotsc, x_n$ but some other sequence $y_0, \dotsc, y_s$ as well. So the existence of a real solution is proved in more general circumstances. A {\it sequence} of homogeneous systems is constructed, each of which can be decomposed into linear factors that are in general position. The latter given condition is an algebraic one.

Each of these systems will have finitely many (hence an odd number) of solutions, with any non-real solution paired with its conjugate. But the equations $F_j$ and $F'_j$ (the new one) can be connected by a homotopy to yield a system valid over some algebraic closure $\Lambda$ of $\RR(t)$. By considering the simplicity of solutions (coming from \cite{vd Waerden, Einf\"uhrung}) in this field, a real solution can be pulled back from the finitely many solutions now seen to exist over the (real-closed) field of real Puiseux series. One should consult the article \cite{Behrend} for details.

We perceive formula (B) as arising, in B\'ezout's Theorem (Thm 4) for a sum of multiplicities of solution rays. Such a situation, for a quantity of equations equal to {\it one less} than the  number of homogeneous variables, would be easier to deal with in case each solution had unit multiplicity. This is indeed the case when  the coefficients are generic (algebraically independent over $\QQ$). At least when we are allowed to  operate over $\CC$ or $\RR$ as coefficients, it would  seem that we could {\it nudge} them one by one  into genericity while homing in on the ``specialized'' solution that we seek over $\RR P^n$.

This strategy best fits the approach from  \cite{Lang, Places} which exhibits {\it both} a ``more algebraic'' and ``more topological'' version to finish off the proof of Theorem 5.
The above-mentioned work of F\. A\. Behrend reduces the problem to one of simple solutions, coming from  a classical criterion for simplicity which we will refer to again. Relevant background is described in \cite{vd Waerden, Einf\"uhrung}, section 39.

The following  gives us a result that would  be a sufficient alternative. It comes from the same textbook of van der Waerden, Dover edition (1945) or Springer-Verlag edition (1973). In section 41 we read ``The intersection of an irreducible $d$-dimensional variety of reduced degree $\gamma$, with a quantity $k \leq d$ generic hypersurfaces of degrees $e_1, e_2,\dotsc, e_k$ respectively, has degree $\gamma \prod_{j=1}^k e_j$. Hence in case $k = d$, this variety consists of (this many) points''.

A similar statement from the earlier book of \cite{Macaulay}, p\. 16 indicates that ``the number of solutions is either $L = l_1 \cdot l_2\cdots l_n$, or infinite, the latter being the case when $F_0$ (a resultant of the system with respect to $x_1, \dotsc, x_n$) vanishes identically''.

Finally, in \cite{Cox AG}, it is proved using an explicit construction that ``the equations $F_1 = \cdots = F_n = 0$ when generic, have $d_1\cdots d_n$ {\it distinct} solutions''. The discussion is in Chapter 3, Section 5, including Exercise 6. The proof involves projective elimination theory and the use of  Macaulay's resultant (which is effective if inefficient). In the sequel it will be seen that we do not need the hard ``generic'' precondition on coefficients to  first finish Theorem 5 and hence the B-U Theorem over $\RR$. We will wish however to avoid the ``Ausnahmefall'' (infinitely many solution-ways). With this in mind, we do use Resultant applications from both \cite{Cox AG} and \cite{vd Waerden, Algebra II}.

We now point out that these transcendental constructions, say in Lang's method can be gotten around in a sense. With the concern that the solutions not be infinite in number, and actually all possess {\it unit} multiplicity, it comes down to whether certain resultants (integer multi-nomials in the coefficients) can possibly vanish. But for given degrees $d_1,\dotsc,d_n$, the ``size'' of these resultants is definitely bounded. Thus we don't need transcendental numbers,  we merely construct sets that are ``sufficiently'' independent. For example, to {\it approximate} $\alpha \in \RR$, we could use $\alpha + \epsilon$ where $\epsilon$ is a small transcendent, on we could use $\epsilon' = (p)^{\frac{1}{q}}$ for large enough primes $p, q \in \NN$. The proof of any of these assertions goes far beyond our intentions.

The alternative offered by Lang at the  end of the 1953 {\it Annals} paper is to use the more familiar mathematics of the standard $\RR$.

In finding a real solution to $F_1(x_0,\dotsc, x_n) = 0$, $F_2(x_0,\dotsc, x_n) = 0$, $F_n(x_0,\dotsc,$\linebreak
 $x_n) = 0$, we have noted several ``algebraic'' proofs of the past, including  those of Macaulay, Behrens, the theory of ``faithful specializations'' (with which van der Waerden replaced a heavy reliance on the explicit use of classical  resultants), Pfister's module-theoretic approach, and finally (in our narrative), the method of real places introduced by S. Lang. We saw how to gain an advantage (through the full complement of simple solutions) by approximating the given coefficients of  $\{F_j\}$ by a set  of algebraically independent coefficients. As \cite{Lang, Places} points out, thus can be done by embedding the real-closed coefficient field $R$ into a real-closed domain $\Omega$ having  many transcendental elements that are infinitesimal with respect to $R$. Such constructions are algebraic and do not use the order-topology of $\RR$.

After mentioning the work of these authors, we assure the loyal Reader that we quickly finish up this  approach to Theorem 5. Multiplicity of solutions is allowed (the coefficients can be specialized), so the  remaining component is an algebraic description of ``multiplicity'' from \cite{vd Waerden, Algebra II}. This is based on the theory of \cite{Kapferer}  and the ``u-resultant''. F\. S\. Macaulay attributes the u-construction to Liouville. Solution of systems by means of variations on the u-resultant figure importantly in Computational Algebra \cite{Cox AG}, \cite{CanMan}, \cite{D'Andrea}.

\bigskip
\head The Real Solution-Ray\endhead

We work with the system of homogeneous equations in $x_0, \dotsc, x_n$ over the standard reals $\RR$
$$F_1 = 0,\dotsc, F_n = 0,$$
although we emphasize the algebraic aspects of  the problem. We avoid the full power of Lang's real-closed domain $\Omega$ by allowing for the metric closeness of $\RR$. We avoid the need to work with systems where each solution has to be simple, and thereby also avoid explicit resultant constructions coming from Elimination Theory.

We do wish to use the methods leading  to the statement of B\'ezout's theorem on page 16 of \cite{vd Waerden, Algebra II}, Section 83. Thus we must ensure that the system (S) possesses only finitely many solutions. We saw how this would come about in case the collection of all coefficients were generic, as it is taken in \cite{Lang, Places}, see also \cite{Cox AG}, Chapter 3. The number of distinct monomials is something like (see \cite{Ryser}),
$$\sum_{j=1}^n \dfrac{(n+d_j)!}{n! d_j!}.$$
\indent
Instead we propose to take {\it all} of these coefficients to lie in $\QQ$ (the rationals), except for one coefficient, which is chosen to be transcendental. Even better, this final real number can be chosen as algebraic but of such an unreachable algebraic order (such as we noted, some $(p)^{\frac{1}{q}}$) that it could never be canceled in the resultant evaluation that arises.

More specially, we may examine ($S$) for ``points at infinity'' by specializing $x_0 = 0$. Now we obtain a system (still homogeneous)
$$\align
 \overline{F}_1(x_1,\dotsc, x_n) &= F_1(0, x_1, \dotsc, x_n) = 0 \\
 \vdots \tag{$\overline{S}$} \\
 \overline{F}_n(x_1,\dotsc, x_n) &= F_n(0, x_1, \dotsc, x_n) = 0
 \endalign
  $$
  in which the number of variables equals the number of equations.

  Hence ($\overline{S}$) is amenable to the theory of Inertial Forms of H. Kapferer (1927).  In our case of interest this boils down to saying that there exists a multi-nomial $R(u_{11}, \dotsc, u_{1n}, \dotsc, u_{nn})$ in the coefficients of  $\overline{S}$ that vanishes precisely when a solution-ray to $\overline{S}$ exists (projective solution). Amongst other properties, $R(u)$ is homogeneous in  the vector of coefficients for $\overline{F}_1$, of total degree $d_2 \cdots d_{n-1}\cdot d_n$, and for $\overline{F}_j$, of total degree $d_1 \cdots \hat{d}_j \cdots d_n$.

  For $R$ to equal $0$ for a particular specialization cannot happen for the case we have chosen of ``all coefficients rational'' (except for the one of them which is chosen transcendental). Hence by this theorem of \cite{Macaulay}, there are no common solutions for $\overline{S}$, hence no solutions ``at infinity'' for $(S)$. A modern and algorithmic account of the Macaulay resultant is available in \cite{Kalorkoti}; see also \cite{Canny}, \cite{CLO} and \cite{Jou}.

  An ideal-theoretic definition and description of solution-multiplicity in given in Chapter XI of \cite{vd Waerden, Algebra II} and in \cite{vd Waerden 1927}. We add to $S$ the linear equation with  independent coefficients
  $$F_0(u) = u_0x_0 + \cdots + u_n\,x_n$$
in order to form the ``$u$-resultant'' of system $S$. The Kapferer (or Inertial) resultant ideal is generated by multi-nomials $b_1(u),\dotsc, b_n(u)$, so that this $b$-system vanishes at $(u_0,\dotsc,u_n)$ exactly when a solution $x = (\xi_0,\dotsc, \xi_n)$ of $S$ exists such that also
  $$L = u_0\, \xi_0 + \cdots + u_n\,\xi_n = 0.$$
  For each such solution $\xi^p = \left(\xi_0^p,\dotsc, \xi_n^p\right)$, we have a linear form $L^p$. Each $b_i(u)$ has roots in the variety defined by $\Lambda = \prod_p L^p (u)$. Since we operate over an algebraically closed field (an extension of $\RR$), we may apply the strong form of Hilbert's Nullstellensatz to obtain $b_i(u)^{\tau_i} \in \Lambda$. Actually the roots of the $b$-system and of $\Lambda  (u)$ are the same so we also have
  $$\Lambda (u)^{\tau} \in \left(b_1(u),\dotsc, b_r(u)\right).$$
  By the theory of Inertial ideals, the  greatest common divisor $R(u) = \mbox{gcd}\left(b_1(u),\right.\dotsc, $\linebreak
  $\left. b_r(u)\right)$ decomposes into the linear factors as indicated:
  $$R(u) = \prod_p L_p^{s_p}(u).$$
  Thus, the linear forms $L^p$ which determine the solution rays of $(S)$ constitute the irreducible factors of the $u$-resultant $R(u)$. The exponents $\{s_p\}$ in the factorization are the solution {\it multiplicities}. Since it is known that the generator $R(u)$ of the (principal) Inertial ideal has total homogeneous degree $D = \prod_{j=1}^n d_j$, we again have B\'ezout's theorem, valid for when the solution-rays for $(S)$ are finite in number:
  $$\sum s_p = D.$$
  \indent
  Now we finish our intended proof that a real solution of $(S)$ exists. The point is that for any {\it non-real} solution $p$, its multiplicity and that of its complex conjugate solution are the same:
  $$\mbox{mult}(p) = \mbox{mult}(\overline{p}),$$
  or $s_p = s_{\overline{p}}$. This observation is made by sheer logic, as the algebraic operations used in calculating $s_p$ do not depend on how an imaginary coordinate was named, $\imath$ or $-\,\imath$. In other words, one may re-label a value $\xi = (-i, \pi + 7i, 4)$ as $\xi' = (i, \pi-7i, 4)$ without affecting the solution algorithm. All ideals, resultants and multiplicities come up again with a superficial change of symbolism. This same fact can be expressed more geometrically of course, as in Chapter IV, 2.2 of \cite{Shafarevich}.

  What remains as far as the use of B\'ezout's theorem is concerned in to see how we have avoided the
  ``Ausnahmefall'' of infinitely many zeros. In that case, B\'ezout's theorem holds true and since the degree $D$ is a  product of odds, and the non-real solutions are paired up, we must obtain a real solution. This again is what is needed in our approach to  the Borsuk-Ulam Theorem and Invariance of Domain.

  The issue of solutions at $\infty$ (where $x_0 = 0$) comes down to the Macaulay resultant taking on a (scalar) value of zero. For the genericity that  we have built into the coefficients of $\overline{S})$, this  is not possible. We picked one coefficient to be transcendental in $\RR$ and the rest algebraic over $\QQ$ (or even rational). Since the resultant construction treats coefficients without prejudice, an equation $\overline{R} (c_{ij}) = 0$ would lead to an algebraic relation not leaving out the chosen ``generic'' one.

  Therefore, given that $(S)$ has no solution-rays at infinity, we infer that the quantity of solution rays is finite. This is a well-known proposition in projective geometry, to which there are several approaches, the more analytical and the more algebraic.

  \bigskip
\head Closed Variety Away from $\infty$\endhead

By making the system $(S)$, $\{F_1,\dotsc, F_n\}$ {\it generic enough}, we avoided solutions (rays) at infinity, so in fact $(S)$ and its associated variety $V$ can be expressed by:
$$\align
 G_1(x_1,\dotsc, x_n) &= F_1(1, x_1, \dotsc, x_n) = 0 \\
  \vdots \tag{$S_1$} \\
  G_n(x_1,\dotsc, x_n) &= F_n(1, x_1, \dotsc, x_n) = 0.
 \endalign
  $$
Hence $V$ is an {\it affine} variety, in particular the $\{G_j\}$ are generally non-homogeneous multi-nomials. We are working in a situation where we need not be concerned with ``real'' fields. The field $K$ of coefficients of $\{F_1,\dotsc, F_n\}$ should be algebraically closed.

We will prove what is required to complete the argument for Theorem 5. The case of interest is where $K = \CC$, so we begin with an argument that uses the order-topology of $\CC$. Subsequently we review an argument from elementary algebraic geometry showing that for any $K$, it is also true that the system $(S_1): \{G_1=0, \dotsc, G_n = 0\}$ also has only finitely many solutions.

Considering first the complex case $K = \CC$, we note that the variety $V$ is a ``projective algebraic set'' and hence {\it compact} in the $\CC$-topology. Now we change the affine coordinates of $\{G_j\}$ if necessary. For parameters $\lambda_j \in \CC,\  j = 1,\dotsc, n-1$ set $x_i' = x_i + \lambda_i x_n$ and $x_n' = x_n$.

\bigskip
\noindent
{\bf Proposition A} \quad A compact, complex affine variety must be a {\it finite set}.

\bigskip
\noindent
{\it Proof}\quad The $\CC(x_1,\dotsc, x_n)$-ideal generated by $\{G_1,\dotsc, G_n\}$ is called $\Cal I (G)$ and its zero-set $\subset \CC^n$ is called $Z(G)$. It is known how to define the ``first elimination ideal'' $\Cal J = \Cal I \cap \CC [x_2, \dotsc, x_n]$.  An induction hypothesis is that ``$Z (\Cal J)$ is {\it bounded} in the $\CC$-norm, implies that $Z (\Cal J)$ is finite''. The base of induction, with one variable say $x_n$, provides of course finitely many solutions.

Now let $g \in \Cal I$ be any multinomial of the ideal. By breaking $g$ up into its homogeneous pieces, it is possible to find parameters $\lambda_1,\dotsc, \lambda_{n-1}$ so that
$$g(x_1',\dotsc, x_n') = \gamma x_1^{'m} + \mbox{lower degree terms in $x_1'$} \tag{\dag}$$
with coefficients $h_\alpha(x_2',\dotsc, x_n')$, where $m$ is the highest total degree of a monomial in $g$.

Such a coefficient $\gamma$ is actually equal to $g_m(1, \lambda_2,\dotsc, \lambda_n)$ where $g_m$ is the homogeneous part of highest degree. In an infinite field, this  expression cannot always equal $0$ unless $g_m$ is  identically $0$, which gives a contradiction.

Next we convert all expressions of the problem into the coordinates $\{x_1',\dotsc, x_n'\}$, then remove the ``prime'' for legibility. At this point we have performed a version of the Noether normalization lemma, see \cite{Arrondo}. Now a ``long solution'' $(c_1, \cdots, c_n) \in Z$ projects into a ``short solution'' $(c_2, \dotsc, c_n) \in Z(\Cal J)$. We could assume that $\Cal J$ gives rise to  an unbounded set of short solutions, or else only a finite quantity of them. If they tend out to infinity, so do the long solutions arising from $(\dag)$. If they are finite in number, $(\dag)$ shows also that the long solutions are finite  in number. Hence if $Z(\Cal J)$ is compact it is finite.\hfill $\blacksquare$

\bigskip

One may phrase this result to say that a ``variety'' over $K$ can only be both affine and projective, when it consists of finitely many solution points. An algebraic set coming from a finitely generated ideal is the union of {\it irreducible} algebraic sets, also called ``varieties'' by some authors. So we may consider a variety $X$ that is also an affine set in $K^n$.

Consider now the field of regular functions
on $X$ consisting of quotients $h/g$ of homogeneous terms $h, g \in K [x_0, x_1,\dotsc, x_n]$ having the same total degree. But $g$ should be non-zero everywhere, so must be a constant, hence also $h$ has to be a constant.

\bigskip
\noindent
{\bf Proposition B} \quad The field of regular functions on an (irreducible) projective variety $X$ is a field of constants $\simeq K$. See \cite{Shafarevich} p. 59.
\bigskip
\noindent
{\it Proof}\quad   Elaborating on our previous argument, we know that $X$ is ``affine'' and hence its ``coordinate ring'' is
$$\Cal O(X) = K[x_1, \dotsc, x_n]/{\Cal I(G)}.$$
But this quotient gives the field $K$ only if  $\Cal I$ is a maximal ideal, which by Hilbert's Nullstellensatz only holds true (over algebraically closed $K$) when $\Cal I$ is the ideal $(x_1-b_1, x_2-b_2,\dotsc, x_n-b_n)$ whose solution zero is the single point $\vec b$, as was to be proved. See \cite{Atiyah}.\hfill $\blacksquare$

\bigskip
Finally we re-work this last result that $X$ must be a finite solution-set, in somewhat greater detail where we employ a ``compactness'' argument modified  from the case of ground field $= \CC$. The new argument applies also to  general (closed) fields. Similar material may be found in a classical exposition, \cite{Shafarevich}.

Consider a regular mapping $f: X \to Y$ of one closed projective set to another. Thus locally, $f$ is defined by a polynomial map.
The {\it graph} of $f$ is the set of pairs $\Gamma_f = \{(x, f(x))\} \subset X \times Y$.

\bigskip
\noindent
{\bf Proposition 1} \quad For a regular mapping $f$, the graph $\Gamma_f$ is (Zariski-) closed in $X \times Y$.

\bigskip
\noindent
{\it Proof}\quad if $\imath$ is the identity $\imath: Y \to Y$, it is seen that $\Gamma_f$ equals the inverse image of $\Gamma_{\imath}$ under $(f, \imath): X \times Y \to X \times Y$,  hence is closed if we know that $\Gamma_{\imath}$ is closed.
But the ``diagonal'' $\Gamma_{\imath} \subset Y \times Y$ is defined by polynomial equations, hence is closed.\hfill $\blacksquare$

\bigskip
\noindent
{\bf Proposition 2} \quad If $X$ is a projective variety, and $Y$ is a projective or affine variety, then the projection $\pi: X \times Y \to Y$ onto the second factor maps closed sets to closed sets.

\bigskip
\noindent
{\bf Remarks} \quad This ``Main Theorem of Elimination Theory'' is covered in textbooks as well as the computational manual \cite{CLO}, Chapter 8, Section 5.

We have referred previously to polynomial conditions (the resultant systems) whose zero-sets define the parameter values (in $Y$) for which a set of equations have solutions in $X$. We saw the following result earlier on.

\bigskip
\noindent
{\bf Corollary 1} \quad If $\varphi$ is a regular function on an {\it irreducible} projective variety then $\varphi  x = c$ for all $x \in X$, so $\varphi$ may be considered as a field element (scalar constant).

\bigskip
\noindent
{\it Proof}\quad Similar to before, $\varphi$ can be viewed as a map  to $\PP^1$ that misses the infinity point. We have from the Proposition that $\varphi(x)$  is closed in $\PP^1$, since $\varphi(x)$ equals the projection to  $\PP^1$ of the graph $\Gamma_{\varphi} \subset X \times \PP^1$. But a closed set in $\AA^1 \subset \PP^1$ is a finite set, which must be a singleton since $X$ is irreducible.\hfill $\blacksquare$

\bigskip
\noindent
{\bf Corollary 2} \quad If a projective set variety $X$ is embedded in an affine $Y$, $X$ consists of finitely many points.

\bigskip
\noindent
{\it Proof}\quad If $Y \subset \AA^m$, the coordinates of image-points of each irreducible component must be constant by Corollary 1. Since there are finitely many components, $X \subset Y$ is a finite set.\hfill $\blacksquare$

\bigskip
This settles again the issue needed for B\'ezout's theorem, that a projective variety avoiding points at infinity must be finite (and $0$-dimensional).

We address this question one final time, letting the Reader pursue the matter further. A {\it finite mapping} $\varphi: X \to Y$ is a regular mapping whose image is (Zariski-) open, and which satisfies an {\it integrality} condition on the induced inclusion of coordinate ring $K[Y] \subset K[X]$. For our purposes, it is enough to know that when $K = \CC$, $\varphi$ must be a finite-to-one continuous mapping of spaces in the $\CC$-topology (a finite covering with some branch points). This is itself a formulation of the Noether Normalization Theorem on $\CC$. First the result:

\bigskip
\noindent
{\bf Proposition 3} \quad Over general $K$, algebraically closed, an irreducible {\it affine} variety $X$ can be mapped to some affine $\AA^m$ by a finite mapping.

\bigskip
\noindent
{\it Proof}\quad See \cite{Shafarevich} p. 65.\hfill $\blacksquare$

\bigskip
We complete our remarks concerning the complex case. For general affine $X$, we use the finite mapping given by Proposition 3 to construct a particular finite mapping $\varphi: X \to \AA^m$. As we saw, for $K = \CC$ such a mapping is continuous and  finite-to-one. In particular $\varphi$ is proper and onto, so if $X$ were compact in the trascendental topology, $\CC^m$ would be too, which gives a contradiction. \hfill $\blacksquare$

\bigskip
The article \cite{Kalorkoti} gives an effective algorithm precisely in the case of ``no zeros at infinity'' and produces the $u$-resultant and in principle its factors. Thus are derived the finitely many solutions, with multiplicity, to the original system $F_1,F_2,\dotsc, F_n$.

\bigskip
\head The B-U Theorem according to Lusternik and Schnirel'mann\endhead

The authors of \cite{L-S} introduced a natural number $\mbox{cat}(M)$ which for our purposes applies to a compact manifold of finite dimension. It turns out that {\bf cat} is actually an invariant of homotopy type \cite{James}. The paper ``M\'ethodes Topologiques...'' seeks to introduce a \text{sharpening} of Morse's \text{inequalities} \cite{Milnor}, and to study geodesics on a Riemannian manifold.

The usual definition of $\cat(M) = k$ is to say that $M$ can be covered by a quantity $k$ open subsets $\{U_i\}$, each of which is contractible to a point ambiently  within $M$ (the inclusion $\imath_i: U_i \to M$ is nil-homotopic).

An inequality $\cat(\RR P^n) \leq n+1$ follows from general considerations of dimension (see below). The more challenging assertion is that $\cat(\RR P^n) \geq n+1$, which follows from the fact that
$$\RR P^1 \subset \RR P^2 \subset \cdots \subset \RR P^n$$
is a chain of similar subspaces, where each inclusion is homologically non-trivial.

The importance of $\cat(\RR P^n) = n+1$ is seen by

\bigskip
\noindent
{\bf Proposition 4}\quad If this calculation holds true, then in every covering of $\SS^n$ by quantity $n+1$ open sets, one of the sets contains an antipodal pain of points $\{x, -x\}$, $x \in \SS^n$. Thus the
Lusternik-Schnirel'mann theorem, see \cite{Matou{\v s}ek}, would be demonstrated.

\bigskip
\noindent
{\it Proof}\quad Let $q: \SS^n \to \RR P^n$ be the canonical double covering (quotient) mapping.  If $\{U_0, \dotsc, U_n\}$ covers $\SS^n$ with no $U_i$ containing any antipodal pair, then $q(U_1),\dotsc, q(U_n)$ must cover $\RR P^n$. Indeed, if $\xi \in \RR P^n$ is not in their union, then for some $y \in U_0$, we have $q(y) = \xi$. Since $-y \notin U_0$, we get  $-y$ belonging to another $U_j,\, j \ne 0$. However, $q(-y) = \xi$ so $\xi \in q(U_j) \subset \bigcup_{i\ne 0} q(U_i)$ after all.  A nil-homotopy in $\SS^n$ of $U_i \subset \SS^n$ induces a nil-homotopy of $q(U_i) \subset \RR P^n$, so we see that $\RR P^n$ has a nil-homotopic cover of size $n$ which gives a contradiction to the hypothesis.\hfill $\blacksquare$

\bigskip

The upper bound we need on $\cat(M)$, that is, one plus the dimension of $M$, is obtained by means of finding a categorical sequence \cite{Fox} for $M$. When $M$ is a finite simplicial complex, it is not difficult to produce such a sequence by means of a ``Balls, Beams, Plates'' construction similar to that used with Haken manifolds. We illustrate this in the specific case where $M $ is a 3-dimensional pseudo-manifold (each two-simplex is the boundary of exactly two three-simplexes). Assume that $M$ is topologically connected.

\bigskip
\centerline{\epsfbox{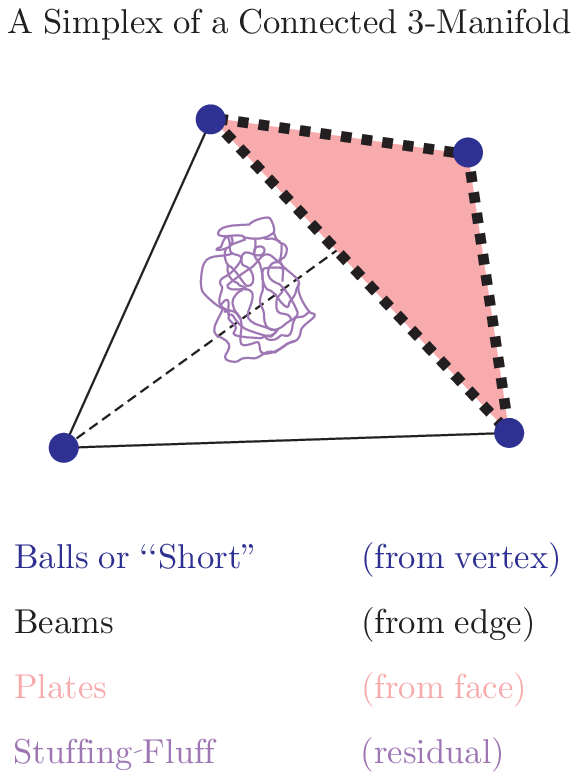}}

\medskip

\centerline{\bf Figure 4}

\bigskip

The vertices of the 3-complex $\Delta$ with  $|\Delta| \simeq M$ are thickened into 3-balls, called ``Shot''. Since $M$ is connected, the Shot is contractible into a point of $M$. Next the 1-simplices are thickened into
``Beams'' which are separated near the vertices, so the collection of Beams is also ambiently contractible. The same holds for the thickened faces, or ``Plates''. Finally, all the open interiors of the 3-simplices of $\Delta$ are united to form the ``Stuffing'' whose inclusion into $M$ is nil-homotopic.

We have covered $M$ with {\it four} contractible open sets, confirming the formula $\cat(M) \leq \dim M+1$. The same method applies to any connected pseudo-manifold of a higher dimension. See Figure 4.

A ``homology'' version $H\cat (M)$ was introduced by \cite{Schnirel'mann}. This number is not greater than       $\cat (M)$. Suppose that for $j \leq 0 , \cdots , n$, $L_j$ is a manifold of dimension $j$, satisfying
$$L_0 \subset L_1 \subset \cdots L_n = M \, ,$$
with the following homological condition: any 1-cycle $\mod \ZZ_2$ of $L_j$ that bounds a $\ZZ_2$-chain in $M$ already bounds in $L_j$. Then we may deduce the following.

\bigskip
\noindent
{\bf Proposition 5}\quad Under the above conditions, it follows that $H \cat (M) \geq n+1$ and hence $\cat (M) \geq n+1$.
\bigskip
\noindent
{\it Proof}\quad See \cite {Fox} and \cite {Schnirel'mann}.\hfill $\blacksquare$

\bigskip
In the case of $M = \RR P^n$, we may define $L_j = \RR P^j$, canonically embedded in $M$, and verify the hypotheses of Proposition 5. Thus we give witness to an earlier proof of the Lusternik-Schnirel'mann theorem, \cite{L-S}, and hence the Borsuk-Ulam theorem, this time based to an extent on ``chain-level intersection'' in homology.

Work in the cohomology ring has largely replaced a historical fashion for chain-level intersection. One shows that the {\it nilpotency index} of the ring gives a lower bound for $\cat (M)$. This integer $k$ is the least such that all $k$-fold cup-products vanish, see \cite{James}. The computation of the ring $H^*(\RR P^n, \ZZ_2)$ seems more involved than the proofs of \cite{Fox} or \cite{Schnirel'mann} sketched above.

In \cite{Goresky-MacPherson}, the authors encourage a return to geometric intersection products as an alternative to the cohomology ring. Another article, by McClure, asserts that these theories are ``probably'' the same as given in the manual \cite{Lefschetz}. Prof. Lefschetz' intersection calculus utilizing ``looping coefficients'' has not often been applied, though there is a monograph \cite{Keller} from  Leipzig (1969) that thoroughly addresses such issues. One hopes that some of the contemporary authorities have read this work. In any case, the old geometric intersection theory (of chains) seems not yet to be fully integrated with a modern homological version. Further discussion can be found in \cite{McClure}.

We suggest that a reworking of the ``intersection-level'' Borsuk-Ulam proof  of \cite{L-S}, based on a specific triangulation and dual triangulation of the real projective space would be of interest, especially to combinatorial mathematicians.

\bigskip
\head Application to Banach Geometry\endhead

The concept of {\it defect} or gap between  two operators on a (real) Banach space, as developed by M\. A\.
Krasnoselskii and co-workers in \cite{KKM}, proved to have fundamental implications concerning the geometry of a Banach space. If $M$ and $N$ are subspaces of finite dimension in a Hilbert space $H$, and $\dim M < \dim N$, then there is  a vector $u \in N$ that is orthogonal to all of $M$. This fact is not hard to see, since in a Hilbert space one can project $M$ into $N$ by a projection $\pi$, where the image is then a linear space of lesser dimension. Some vector $u \in N$ that is orthogonal to $\mbox{Im}(\pi(M))$ will then also be orthogonal to $M$ itself.

If alternatively $M$ and $N$ are subspaces of a normed linear (or Banach) space, the analogous result is less obvious. For one thing, it is necessary to define the ``orthogonality''  of a given vector $u$ with some subspace $M$. We may adopt the definition
$$d(u, M) \equiv \inf \left\{\|u-y\|:y\in M\right\}.$$
Thus the distance from $u$ to the subspace should be minimized as the distance to the 0 subspace, giving the norm of $u$, that is, $\|u\|$.

Call this result (the existence of $u \in N$ orthogonal to $M \subset N$) the Theorem on the Deviation of Subspaces \cite{Brown}.
In fact the statement is logically equivalent (by a short derivation) to the Borsuk-Ulam theorem.

We indicate some features of the proof of Deviation of Subspaces from the B-U theorem. Without loss of generality, one may assume that $\dim N = \dim M+1$. For a first case take it that the Banach space $E$ is just the (finite-dimensional) {\it sum} of $M$ and $N$, and that $E$ is {\it strictly convex}.  This means that for two {\it linearly independent} vectors $u, v \in E$, we have $\|u+v\| < \|u\| + \|v\|$. Now a derivation from the elementary theory of normed vector spaces shows that every $u \in E$ has a {\it nearest} vector $\psi (u) \in M$, and that $\psi: E \to M$ is continuous in the norm topology. In case $E$ is not a Hilbert space, $\psi$ might not be a {\it linear} mapping, but it does satisfy
$$\psi(-u) = -\psi(u),$$
so is antipode-preserving on the sphere of ``norm one'' vectors of $N$. Hence, by the Borsuk-Ulam theorem, see \cite{Matou{\v s}ek}, there exists $u \in N$ with $\|u\|=1$ and $\psi (u) = 0$. As alluded to above, this vector in $N$ is the one we seek, it is orthogonal in the Banach sense, to all of $M$.
For the general case where $E$ is not strictly convex, given $\epsilon > 0$, the experts (see \cite{Gohberg-Krein}) construct a new metric $\|\quad\|_0$ on $E$ which satisfies $$\|v\| \leq \|v\|_0 \leq (1+\epsilon)\|v\|$$ for all $v \in E$. It turns out that the sphere $\{v: \|v\|_0 = 1\}$ is strictly convex. With the new norm, one can find $u$ of $\mbox{norm}=1$ that is orthogonal to $M$.  

Actually, this $u$ depends on the choice of norm and should be written $u_{\epsilon}$. As $\epsilon \to 0$, one picks out a convergent subsequence of the $u$,  where $\epsilon= 2^{-k}$ (the ``original'' norms of these vectors go to 1), and this vector is then shown to be orthogonal to $M$.\hfill $\blacksquare$
\bigskip 
A\. L\. Brown proved the converse in \cite{Brown}. We already have a proof of the B-U theorem, but he applies Deviation of Subspaces to the space $E = C(\SS^n)$ of continuous real-valued functions on the $n$-sphere, equipped with the ``supremum'' (or ``uniform'') norm. Let $N \subset E$ be generated by the  coordinate functions of $\RR^{n+1}$, $\SS^n \subset \RR^{n+1}$. Let $M$ be generated by the $n$ coordinate functions on $\RR^n$, after applying $\varphi: \SS^n \to \RR^n$, a continuous, antipode-preserving mapping. 

 One need only show that there is a vector $w \in \SS^n$ with $\varphi (w) = \vec{0} \in \RR^n$. But if $z \in N$ can be found, orthogonal to $M$ as asserted by the Deviation theorem, $z$ is actually a {\it linear functional} on $\RR^{n+1}$ that attains its norm in $C(\SS^n)$ at a (unique) antipodal pair $\{w , -w\}$. This choice of $w \in \SS^n$ turns out to provide the ``Borsuk-Ulam'' vector that is required.\hfill $\blacksquare$

\newpage
\end{document}

\end{document}